\documentclass[11pt,a4paper]{amsart}
\begin{document}
\title{Hypergeometric identities for 10 extended Ramanujan-type series}
\author{Jes\'{u}s Guillera}
\address{Zaragoza, (Spain)}
\email{jguillera@gmail.com}
\date{}

\maketitle

\begin{abstract}
We prove, by the WZ-method, some hypergeometric identities which relate ten extended Ramanujan type series to simpler
hypergeometric series. The identities we are going to prove are valid for all the values of a parameter $a$ when they are convergent. Sometimes, even if they do not converge, they are valid if we consider these identities as limits.
\end{abstract}

\section{Introduction}
In this paper we are going to prove some hypergeometric identities by using WZ-pairs [10], [11], that is,
functions $F(n,k)$ and $G(n,k)$ related by
\[ G(n,k+1)-G(n,k)=F(n+1,k)-F(n,k). \]
A package, written by D. Zeilberger and called {\it ekhad} [8],
allows one to obtain $F$ from $G$ or viceversa. For our purposes we need the following implication:

\begin{align}
G(n,k+1)-G(n,k) &=F(n+1,k)-F(n,k) \nonumber \\
\Longrightarrow G(n+a,k+1)-G(n+a,k) &=F(n+a+1,k)-F(n+a,k). \nonumber
\end{align}
So, if we denote $F_a(n,k)=F(n+a,k)$ and $G_a(n,k)=G_(n+a,k)$ then obviously $F_a(n,k)$ and $G_a(n,k)$ is an WZ-pair for every value of $a$. So, we can write:
\[ G_a(n,k+1)-G_a(n,k)=F_a(n+1,k)-F_a(n,k) \]
and summing from $n=0$ to $\infty$ we get
\[ \sum_{n=0}^{\infty}[G_a(n,k+1)-G_a(n,k)]=\sum_{n=0}^{\infty}[F_a(n+1,k)-F_a(n,k)]=-F_a(0,k) \]
which implies that
\begin{align}
\sum_{n=0}^{\infty} G_a(n ,0)
&= \sum_{n=0}^{\infty} G_a(n,1)+ F_a(0,0)=\sum_{n=0}^{\infty} G_a(n,2)+ F_a(0,1)+F_a(0,0) \nonumber \\
&= \sum_{n=0}^{\infty} G_a(n,3)+ F_a(0,2)+F_a(0,1)+F_a(0,0) \nonumber \\
&= \sum_{n=0}^{\infty} G_a(n,4)+\sum_{k=0}^{3} F_a(0,k) \nonumber
\end{align}
and continuing the recursion we arrive to [2]:
\[ \sum_{n=0}^{\infty} G_a(n ,0)= \lim_{k \to \infty} \sum_{n=0}^{\infty} G_a(n,k)+ \sum_{k=0}^{\infty} F_a(0,k), \]
which is the result we are going to use to get the identities. Finally some easy transformations of $G_a(n,0)$ and $F_a(0,k)$ into rising factorials lead to the desire form of the identities and the proof is complete.

\section{Identities for some extended Ramanujan series}
We consider an extension with a variable $a$ of some Ramanujan series for $1/\pi$ [3], [4] and [9] and get some identities. Once we find an adequate WZ-pair then the difficulty of the proof consist of getting the function
\[ S(a)=\lim_{k \to \infty} \sum_{n=0}^{\infty} G_a(n,k) \]
and very often we have guessed it from
\[ S(a)=\sum_{n=0}^{\infty} G_a(n ,0)-\sum_{k=0}^{\infty} F_a(0,k) \] without finding the proof, so our proofs are incomplete. We will
explain with more detail identities (1) and (2). For the other identities we will limit ourselves to indicate the WZ-pairs which essentially encapsulate each proof because they are all similar.

\subsection{Identity 1}
Let
\[ f(a)=\sum_{n=0}^{\infty} {1 \over 2^{2n}} \frac{ \left( a +\frac{1}{2}\right)_n^3}{(a+1)_n^3}[6(n+a)+1] \]
then we have
\begin{equation}\label{}
f(a)=8a \sum_{n=0}^{\infty} \frac{ \left( \frac{1}{2} \right)_n^2}{(a+1)_n^2}
\end{equation}
and also
\begin{equation}\label{}
f(a)={4 \over \pi} \cdot  { 4^a \over \cos^2  \pi a } \cdot {1_a^3\over \left( {1 \over 2} \right)_a^3}  +  {16 a^2 \over 2a-1}\sum_{n=0}^{\infty}  \frac{  \left( \frac{1}{2} \right)_n  \left(a+\frac{1}{2} \right)_n} {(a+1)_n \left( {3 \over 2} -a \right)_n}.
\end{equation}
From (1) we easily get the evaluation:
\[ f \left(\frac{1}{2}\right)={\pi^2 \over 2} \]
and from (2)
\[
f(0)={4 \over \pi}, \quad \quad \quad f'(0)={32  \over \pi} \ln 2, \quad \quad \quad f''(0)=
{4 \over \pi} \left( 64 \ln^2 2 -3 \pi^2 \right).
\]
\begin{proof}
The proof of (1) is essentially encapsulated in the following WZ-pair:
\[ F(n,k)=8 \cdot B(n,k) \cdot n, \]
\[ G(n,k)=B(n,k) \cdot (6n+4k+1), \]
where
\[ B(n,k)=\frac{1}{2^{8n} \cdot 2^{4k}} \cdot \frac{(2k)!^2 \cdot (2n)!^3}{(n+k)!^2 \cdot k!^2 \cdot n!^4}. \]
For proving it we get:
\[
\sum_{n=0}^{\infty} G_a(n,0)=\sum_{n=0}^{\infty} G(n+a,0)={ \left( {1 \over 2} \right)^3_a \over 1_a^3} \cdot
{1 \over 4^a} \cdot \sum_{n=0}^{\infty} {1 \over 2^{2n}} \frac{\left(a+\frac{1}{2}\right)_n^3}{(a+1)_n^3}[6(n+a)+1],
\]
\[
\sum_{n=0}^{\infty} F_a(0,k)=\sum_{n=0}^{\infty} F(a,k)={ \left( {1 \over 2} \right)^3_a \over 1_a^3} \cdot
{1 \over 4^a} \cdot \sum_{n=0}^{\infty} 8a \cdot \frac{ \left( \frac{1}{2} \right)_n^2}{(a+1)_n^2}.
\]
If $a=0$ then
\[ \lim_{k \to \infty} \sum_{n=0}^{\infty} G(n,k)=\lim_{k \to \infty} G(0,k)={4 \over \pi} \]
and if $a>0$ then
\[ \lim_{k \to \infty} \sum_{n=0}^{\infty} G_a(n,k)=0 \]
and the proof of (2) is essentially encapsulated in the WZ-pair:
\[ F(n,k)=B(n,k) \cdot {16 \cdot n^2 \over 2n-2k-1}, \]
\[ G(n,k)=B(n,k) \cdot (6n+2k+1), \]
where
\[
B(n,k)={(-1)^k \over 2^{8n} \cdot 2^{4k}} \cdot {(2k)! \cdot (2n+2k)! \cdot \left(n-k-{1 \over 2} \right)! \cdot
(2n)!^2 \over k! \cdot (n+k)!^2 \cdot \left(n-{1 \over 2} \right)! \cdot n!^4}.
\]
For proving it we get:
\[
\sum_{n=0}^{\infty} G_a(n,0)=\sum_{n=0}^{\infty} G(n+a,0)={ \left( {1 \over 2} \right)^3_a \over 1_a^3} \cdot {1 \over 4^a} \cdot \sum_{n=0}^{\infty} {1 \over 2^{2n}} \frac{\left(a+\frac{1}{2}\right)_n^3}{(a+1)_n^3}[6(n+a)+1]
\]
\[
\sum_{n=0}^{\infty} F_a(0,k)=\sum_{n=0}^{\infty} F(a,k)={ \left( {1 \over 2} \right)^3_a \over 1_a^3} \cdot
{1 \over 4^a} \cdot {16 a^2 \over 2a-1} \cdot \sum_{k=0}^{\infty} \frac{  \left( \frac{1}{2} \right)_k  \left( a+\frac{1}{2} \right)_k} {(a+1)_k \left( {3 \over 2} -a \right)_k}
\]
and we have guessed that
\[ S(a)=\sum_{n=0}^{\infty} G_a(n,0)-\sum_{k=0}^{\infty} F_a(0,k)={4 \over \pi} \cdot {1 \over \cos^2 \pi a }. \]
\end{proof}

\subsection{Identity 2}
Let
\[ f(a)=\sum_{n=0}^{\infty} {1 \over 2^{6n}} \frac{ \left( a +\frac{1}{2}\right)_n^3}{(a+1)_n^3}[42(n+a)+5], \]
then we have
\begin{equation}\label{}
f(a)=32a \sum_{n=0}^{\infty} \frac{ \left( a+\frac{1}{2} \right)_n^2}{(2a+1)_n^2}
\end{equation}
and also
\begin{equation}\label{}
f(a)={16 \over \pi} \cdot  { 64^a \over \cos^2  \pi a } \cdot {1_a^3 \over \left( {1 \over 2} \right)_a^3}  +  {128 a^2 \over 2a-1} \sum_{n=0}^{\infty}  \frac{  \left( a+\frac{1}{2} \right)_n^2} {(2a+1)_n \left( {3 \over 2} -a \right)_n}.
\end{equation}
From (3) we easily get the evaluation:
\[ f \left(\frac{1}{2}\right)={8 \pi^2 \over 3} \]
and from (4)
\[
f(0)={16 \over \pi}, \quad \quad \quad f'(0)={192  \over \pi} \ln 2, \quad \quad \quad
f''(0)={16 \over \pi} \left(144 \ln^2 2  -7 \pi^2 \right).
\]
\begin{proof}
The proof of (3) is essentially encapsulated in the following WZ-pair:
\[ F(n,k)=32 \cdot B(n,k) \cdot n, \]
\[ G(n,k)=B(n,k) \cdot {(2n+2k+1)^2 \cdot (42n+4k+5)-32 \cdot k \cdot n \cdot (4n+3k+2) \over (2n+k+1)^2}, \]
where
\[
B(n,k)=\frac{1}{2^{12n} \cdot 2^{4k}} \cdot \frac{(2n+2k)!^2 \cdot (2n)!^3}{(2n+k)!^2 \cdot (n+k)!^2 \cdot
n!^4}
\]
and the proof of (4) is essentially encapsulated in the
WZ-pair:
\[ F(n,k)=B(n,k) \cdot {128 \cdot n^2 \over 2n-2k-1}, \]
\[ G(n,k)=B(n,k) \cdot {(2n+2k+1) \cdot (42n+2k+5)-32 \cdot k \cdot n \over 2n+k+1}, \]
where
\[ B(n,k)={(-1)^k \over 2^{12n} \cdot 2^{4k}} \cdot {(2n+2k)!^2 \cdot \left(n-k-{1 \over 2} \right)! \cdot
(2n)!^2 \over (n+k)!^2 \cdot (2n+k)! \cdot \left(n-{1 \over 2}
\right)! \cdot n!^4}.
\]
\end{proof}

\subsection*{Identity 3}
Let
\[ f(a)=\sum_{n=0}^{\infty} {(-1)^n \over 2^{3n}} \frac{\left( a +\frac{1}{2}\right)_n^3}{(a+1)_n^3}[6(n+a)+1], \]
then we have
\begin{equation}\label{}
f(a)=4a \sum_{n=0}^{\infty} \frac{ \left( \frac{a}{2}+\frac{1}{4} \right)_n
\left( \frac{a}{2}+\frac{3}{4} \right)_n}{(a+1)_n^2}
\end{equation}
and also
\begin{equation}\label{}
f(a)={2 \sqrt{2} \over \pi} \cdot  { 8^a \over \cos  \pi a } \cdot {1_a^3 \over \left( {1 \over 2} \right)_a^3}  +  {16 a^2 \over 2a-1} \sum_{n=0}^{\infty}  {1 \over 2^n} \cdot \frac{ \left( a+\frac{1}{2} \right)_n^2} {(a+1)_n \left( {3 \over 2} -a \right)_n}.
\end{equation}
From (5) we easily get the evaluation:
\[ f \left(\frac{1}{2}\right)=4G \]
and from (6)
\[
f(0)={2 \sqrt{2} \over \pi}, \quad \quad \quad f'(0)={18 \sqrt{2}  \over \pi} \ln 2, \quad \quad \quad
f''(0)={2 \sqrt 2 \over \pi} \left(81  \ln^2 2 -4  \pi^2 \right).
\]
\begin{proof}
The proof of (5) is essentially encapsulated in the following WZ-pair:
\[ F(n,k)=4 \cdot B(n,k) \cdot n, \]
\[ G(n,k)=B(n,k) \cdot(6n+4k+1), \]
where
\[ B(n,k)=\frac{(-1)^n}{2^{9n} \cdot 2^{6k}} \cdot \frac{(2n+4k)! \cdot (2n)!^2}{(n+2k)! \cdot (n+k)!^2 \cdot n!^3} \]
and the proof of (6) is essentially encapsulated in the WZ-pair:
\[ F(n,k)=B(n,k) \cdot {16 \cdot n^2 \over 2n-2k-1}, \]
\[ G(n,k)=B(n,k) \cdot (6n+2k+1), \]
where
\[
B(n,k)={(-1)^k \over 2^{9n} \cdot 2^{5k}} \cdot {(2n+2k)!^2 \cdot \left(n-k-{1 \over 2} \right)! \cdot (2n)!
\over (n+k)!^3 \cdot \left(n-{1 \over 2} \right)! \cdot n!^3}.
\]
\end{proof}

\subsection*{Identity 4}
Let
\[
f(a)=\sum_{n=0}^{\infty} {(-1)^n \over 2^{2n}} \frac{ \left( a +\frac{1}{2}\right)_n \left( a +\frac{1}{4}\right)_n
\left( a +\frac{3}{4}\right)_n}{(a+1)_n^3}[20(n+a)+3],
\]
then we have
\begin{equation}\label{}
f(a)=16a \sum_{n=0}^{\infty} \frac{ \left( \frac{1}{2} \right)_n \left( a+\frac{1}{2} \right)_n}{(a+1)_n (2a+1)_n}
\end{equation}
and also
\begin{equation}\label{}
f(a)={8 \over \pi} \cdot  { 4^a \over \cos  \pi a } \cdot {1_a^3 \over \left( {1 \over 2} \right)_a \left( {1 \over 4} \right)_a \left( {3 \over 4} \right)_a} +  {48 a^2 \over 2a-1} \sum_{n=0}^{\infty} {1 \over 4^n} \cdot { \left( a+\frac{1}{2} \right)_n \left( 2a+\frac{1}{2} \right)_n \over {(a+1)_n \left( {3
\over 2} - a \right)_n}}.
\end{equation}
From (7) we easily get the evaluation:
\[ f \left(\frac{1}{2}\right)=16 \cdot \ln{2} \]
and from (8)
\[
f(0)={8 \over \pi}, \quad \quad \quad f'(0)={80  \over \pi}\ln 2, \quad \quad \quad
f''(0)={8 \over \pi} \left(100 \ln^2 2 -5 \pi^2 \right).
\]
\begin{proof}
The proof of (7) is essentially encapsulated in the following
WZ-pair:
\[ F(n,k)=4 \cdot B(n,k) \cdot n, \]
\[ G(n,k)=B(n,k) \cdot {(2n+2k+1) \cdot (20n+4k+3)-16 \cdot k \cdot n \over 2n+k+1}, \]
where
\[
B(n,k)=\frac{(-1)^n}{2^{10n} \cdot 2^{4k}} \cdot \frac{(2k)! \cdot (2n+2k)! \cdot (4n)!} {(2n+k)!
\cdot (n+k)!^2 \cdot k! \cdot n!^2}
\]
and the proof of (8) is essentially encapsulated in the WZ-pair:
\[ F(n,k)=B(n,k) \cdot {48 \cdot n^2 \over 2n-2k-1}, \]
\[ G(n,k)=B(n,k) \cdot {(2n+2k+1) \cdot (20n+2k+3)-24 \cdot k \cdot n \over 2n+1}, \]
where
\[
B(n,k)={(-1)^k \over 2^{10n} \cdot 2^{6k}} \cdot {(2n+2k)! \cdot (4n+2k)! \cdot \left(n-k-{1 \over 2} \right)!
\over (2n+k)! \cdot (n+k)!^2 \cdot \left(n-{1 \over 2} \right)!
\cdot n!^2}.
\]
\end{proof}

\subsection{Identity 5}
Let
\[
f(a)=\sum_{n=0}^{\infty} {(-1)^n 3^{3n}\over 2^{9n}} \frac{\left( a +\frac{1}{2}\right)_n \left( a +\frac{1}{6}\right)_n
\left( a +\frac{5}{6}\right)_n}{(a+1)_n^3}[154(n+a)+15],
\]
then we have
\begin{equation}\label{}
f(a)=128a \sum_{n=0}^{\infty} \frac{ \left(\frac{a}{2}+\frac{1}{4} \right)_n \left( \frac{a}{2}+\frac{3}{4}
\right)_n}{(a+1)_n (2a+1)_n}
\end{equation}
and also
\begin{equation}\label{}
f(a)={32 \sqrt{2} \over \pi} \cdot  { 512^a \over 27^a \cdot \cos \pi a} \cdot {1_a^3 \over \left( {1 \over 2} \right)_a \left( {1 \over 6} \right)_a \left( {5 \over 6} \right)_a} +  {512 a^2 \over 2a-1}
\sum_{n=0}^{\infty} {1 \over 2^n}  { \left( a+\frac{1}{2} \right)_n \left( 3a+\frac{1}{2} \right)_n \over {(2a+1)_n \left(
{3 \over 2} - a \right)_n}}.
\end{equation}
From (9) we easily get the evaluation:
\[ f \left(\frac{1}{2}\right)=128 \ln{2} \]
and from (10)
\[
f(0)={32 \sqrt{2} \over \pi}, \quad \quad \quad f'(0)={480 \sqrt{2}  \over \pi} \ln 2, \quad \quad \quad
f''(0)={32 \sqrt 2 \over \pi} \left(225 \ln^2 2 -11 \pi^2 \right).
\]
\begin{proof}
The proof of (9) is essentially encapsulated in the following
WZ-pair:
\[ F(n,k)=128 \cdot B(n,k) \cdot n, \]
\[ G(n,k)=B(n,k) \cdot {(2n+4k+1) \cdot (154n+16k+15)-384 \cdot k \cdot n \over 2n+k+1}, \]
where
\[
B(n,k)=\frac{(-1)^n} {2^{15n} \cdot 2^{6k}} \cdot \frac{(2n+4k)! \cdot (6n)!} {(2n+k)! \cdot (n+2k)! \cdot (n+k)! \cdot n! \cdot (3n)!}
\]
and the proof of (10) is essentially encapsulated in the WZ-pair:
\[ F(n,k)=B(n,k) \cdot {512 \cdot n^2 \over 2n-2k-1}, \]
\begin{align}
G(n,k)=B(n,k) \cdot \left[ {(2n+2k+1) \cdot (6n+2k+3) \over (2n+1) \cdot (6n+3k+3)} \cdot (154n+6k+15) \right. \nonumber \\ \left.  - {32 kn \over 3} \cdot  {38n+14k+19 \over (2n+1) \cdot (2n+k+1)} \right], \qquad \qquad \qquad \nonumber
\end{align}
where
\[
B(n,k)={(-1)^k \over 2^{15n} \cdot 2^{5k}} \cdot {(2n+2k)! \cdot (6n+2k)! \cdot \left(n-k-{1 \over 2} \right)!
\over (2n+k)! \cdot (3n+k)! \cdot (n+k)! \cdot \left(n-{1 \over 2} \right)! \cdot n!^2}.
\]
\end{proof}

\subsection*{Identity 6}
Let
\[
f(a)=\sum_{n=0}^{\infty} {1 \over 3^{2n}} \frac{ \left( a +\frac{1}{2}\right)_n \left( a +\frac{1}{4}\right)_n
\left( a +\frac{3}{4}\right)_n}{(a+1)_n^3}[8(n+a)+1],
\]
then we have
\begin{equation}\label{}
f(a)={2 \sqrt{3} \over \pi} \cdot  { 9^a \over \cos  2 \pi a } \cdot {1_a^3 \over \left( {1 \over 2} \right)_a \left( {1 \over 4} \right)_a \left( {3 \over 4} \right)_a} +  {36 a^2 \over 4a-1} \sum_{n=0}^{\infty} \left( {3 \over 4} \right)^n \frac{ \left( \frac{1}{2} \right)_n \left( a+\frac{1}{2} \right)_n}{(a+1)_n \left( {3 \over 2} -2a \right)_n}.
\end{equation}
From (11) we easily get the evaluation:
\[ f \left(\frac{1}{2}\right)={\sqrt 3 \cdot \pi} \]
and also from (11) we get
\[
f(0)={2 \sqrt{3}  \over \pi}, \quad \quad \quad f'(0)={4 \sqrt 3 \over \pi} \left( \ln 3+4 \ln 2 \right),
\quad \quad \quad $$ $$f''(0)={4 \sqrt 3 \over \pi} \left( 32
\ln^2 2 + 2 \ln^2 3 +16  \ln 3 \ln 2 - 3 \pi^2 \right).
\]
\begin{proof}
It is essentially encapsulated in the following WZ-pair:
\[ F(n,k)=B(n,k) \cdot {36 \cdot n^2 \over 4n-2k-1}, \]
\[ G(n,k)=B(n,k) \cdot (8n+2k+1), \]
where
\[
B(n,k)={(-1)^k \cdot 3^k \over 2^{8n} \cdot 3^{2n} \cdot 2^{6k}} \cdot {(2k)! \cdot (2n+2k)! \cdot \left(2n-k-{1 \over 2} \right)! \cdot (4n)! \over k! \cdot (n+k)!^2 \cdot \left(2n-{1 \over 2} \right)! \cdot (2n)! \cdot n!^2}.
\]
\end{proof}

\subsection*{Identity 7}
Let
\[
f(a)=\sum_{n=0}^{\infty} {(-1)^n \over 2^{4n} \cdot 3^{n}} \frac{ \left( a +\frac{1}{2}\right)_n \left( a +\frac{1}{4}\right)_n\left( a +\frac{3}{4}\right)_n}{(a+1)_n^3}[28(n+a)+3],
\]
then we have
\begin{equation}\label{}
{16 \sqrt{3} \over 3 \pi} \cdot  { 48^a \over \cos  \pi a } \cdot {1_a^3 \over \left( {1 \over 2} \right)_a \left( {1 \over 4} \right)_a \left( {3 \over 4} \right)_a} +  {96 a^2 \over 2a-1} \sum_{n=0}^{\infty} \left( {3 \over 4} \right)^n   { \left( a+\frac{1}{2} \right)_n  \left( 2a+\frac{1}{2} \right)_n \over {(2a+1)_n \left( {3 \over 2} - a \right)_n}}.
\end{equation}
From (12) we get:
\[
f(0)={16 \sqrt{3}  \over 3 \pi}, \quad \quad \quad f'(0)={16 \sqrt 3 \over 3 \pi} \left(  \ln 3+ 12 \ln 2 \right),
\quad \quad \quad
\]
\[
f''(0)={16 \sqrt 3 \over 3 \pi} \left( 144 \ln^2 2+ \ln^2 3 +24  \ln 3 \ln 2 -9  \pi^2 \right).
\]
\begin{proof}
It is essentially encapsulated in the following WZ-pair:
\[ F(n,k)=B(n,k) \cdot {96 \cdot n^2 \over 2n-2k-1}, \]
\[ G(n,k)=B(n,k) \cdot {(2n+2k+1) \cdot (28n+2k+3)-24 \cdot k \cdot n \over 2n+k+1}, \]
where
\[
B(n,k)={(-1)^k \cdot 3^k \cdot \over 2^{12n} \cdot 3^{n} \cdot 2^{6k}} \cdot {(2n+2k)! \cdot (4n+2k)! \cdot \left(n-k-{1 \over 2} \right)! \cdot (2n)! \over (2n+k)!^2 \cdot (n+k)! \cdot \left(n-{1 \over 2} \right)! \cdot n!^3}.
\]
\end{proof}

\section{Identities for a new kind of extended series}
We consider an extension with a variable $a$ of some new series for $1/\pi^2$ obtained by the author [5-7] and get some
identities. Once we find an adequate WZ-pair then the difficulty of the proof consist of getting the function
\[ S(a)=\lim_{k \to \infty} \sum_{n=0}^{\infty} G_a(n,k) \]
but fortunately for the identities in this section we can commute the limit with the sum allowing us to obtain the function $S(a)$. We will prove with full detail identity (14). For the other identities we will limit ourselves to indicate the WZ-pairs which essentially encapsulate each proof because they are all similar.

\subsection{Identity 8}
Let $f(a)$ and $g(a)$ be the functions:
\[
f(a)=\sum_{n=0}^{\infty} {(-1)^{n} \over 2^{2n}} \frac{ \left( a+\frac{1}{2}\right)_n^5}{(a+1)_n^5}[20(n+a)^2+8(n+a)+1],
\]
\[ g(a)=\sum_{n=0}^{\infty} {1 \over 2^{2n}} \frac{ \left( a +\frac{1}{2}\right)_n^3}{(a+1)_n^3}[6(n+a)+1)], \]
then we have
\begin{equation}\label{}
f(a)=8a \sum_{n=0}^{\infty} \frac{ \left( \frac{1}{2} \right)_n^4}{(a+1)_n^4} (4n+2a+1)
\end{equation}
and
\begin{equation}\label{}
f(a)= {2 \over \pi} \cdot {1 \over \cos \pi a} \cdot { 1_a ^2 \over { \left(1 \over 2 \right)}_ a^2} \cdot g(a) +{32 a^3 \over 2a-1} \sum_{n=0}^{\infty}  \frac{  \left( \frac{1}{2} \right)_n^2 \left( a+\frac{1}{2} \right)_n} {(a+1)_n^2 \left( {3 \over 2} -a \right)_n}
\end{equation}
and also
\begin{align}\label{}
f(a)={8 \over \pi^2} \cdot {4^a \over \cos^3 \pi a} { 1_a ^5 \over { \left(1 \over 2 \right)}_ a^5} &+ {32 \over \pi} \cdot {1 \over \cos  \pi a} \cdot { 1_a ^2 \over { \left(1 \over 2 \right)}_ a^2} \cdot {a^2 \over 2a-1} \cdot \sum_{n=0}^{\infty}  \frac{  \left( \frac{1}{2} \right)_n \left( a+\frac{1}{2} \right)_n} {(a+1)_n \left( {3 \over
2} -a \right)_n} \\
&+{32 a^3 \over 2a-1} \sum_{n=0}^{\infty}  \frac{  \left( \frac{1}{2} \right)_n^2 \left( a+\frac{1}{2} \right)_n} {(a+1)_n^2 \left( {3 \over 2} -a \right)_n}. \nonumber
\end{align}
From (13) we easily get the evaluation:
\[ f \left(\frac{1}{2}\right)=7 \zeta(3). \]
From (14) we get
\[
\lim_{a \to 0} \left[ {f(a) \over a^3} -{2 \over \pi} \cdot {1 \over \cos \pi a} \cdot { 1_a ^2 \over { \left(1 \over 2 \right)}_ a^2} \cdot {g(a) \over a^3} \right] =-32 \cdot \sum_{n=0}^{\infty} { {2n \choose n}^2 \over 16^n \cdot (2n+1)} = -128 \cdot {G \over \pi}.
\]
And from (15) we get
\[
f(0)={8 \over \pi^2}, \quad \quad \quad f'(0)={96 \over \pi^2} \ln 2, \quad \quad \quad
f''(0)={64 \over 3 \pi^2} \left( 54 \ln^2 2 - \pi^2 \right).
\]
\begin{proof}
The proof of (13) is essentially encapsulated in the WZ-pair:
\[ F(n,k)=B(n,k) \cdot 8n \cdot (2n+4k+1), \]
\[ G(n,k)=B(n,k) \cdot (20n^2+8n+1+24kn+8k^2+4k), \]
where
\[ B(n,k)=\frac{(-1)^n}{2^{12n} \cdot 2^{8k}} \cdot \frac{(2k)!^4 \cdot (2n)!^5}{(n+k)!^4 \cdot k!^4 \cdot n!^6}. \]
The proof of (14) is essentially encapsulated in the WZ-pair:
\[ F(n,k)=B(n,k) \cdot {32 \cdot n^3 \over 2n-2k-1}, \]
\[ G(n,k)=B(n,k) \cdot (20n^2+12kn+8n+2k+1), \]
where
\[
B(n,k)={(-1)^k \cdot (-1)^n \over 2^{12n} \cdot 2^{6k}} \cdot {(2n)!^4 \cdot (2n+2k)! \cdot (2k)!^2 \cdot \left(n-k-{1 \over 2} \right)! \over n!^7 \cdot (n+k)!^3 \cdot k!^2 \cdot \left(n-{1 \over 2} \right)!}.
\]
For proving it we get:
\[
\sum_{n=0}^{\infty} G(n+a,0)={\left( {1 \over 2}  \right)_a^5 \over 1_a^5} \cdot {\cos \pi a \over 4^a } \cdot \sum_{n=0}^{\infty} {(-1)^{n} \over 2^{2n}} \frac{ \left( a+\frac{1}{2}\right)_n^5}{(a+1)_n^5}[20(n+a)^2+8(n+a)+1)],
\]
\[
\sum_{k=0}^{\infty} F(a,k)={\left( {1 \over 2}  \right)_a^5 \over 1_a^5} \cdot {\cos \pi a \over 4^a } \cdot {32 a^3 \over 2a-1} \cdot \sum_{k=0}^{\infty} \frac{ \left( \frac{1}{2} \right)_k^2 \left( a+\frac{1}{2} \right)_k} {(a+1)_k^2 \left( {3 \over 2} -a \right)_k}.
\]
For obtaining the limit:
\[ \lim_{k \to \infty} \sum_{n=0}^{\infty} G_a(n,k), \]
we proceed as follows:
\begin{align}
\lim_{k \to \infty} \sum_{n=0}^{\infty} G_a(n,k)
&= \lim_{k \to \infty} \sum_{n=0}^{\infty} G_a(n,k+1) \nonumber \\
&= \lim_{k \to \infty} \sum_{n=0}^{\infty} \left[ G_a(n,0) \prod_{j=0}^{k} {G_a(n,j+1) \over G_a(n,j)} \right] \nonumber \\
&=\sum_{n=0}^{\infty} \lim_{k \to \infty} \left[ G_a(n,0) \prod_{j=0}^{k} {G_a(n,j+1) \over G_a(n,j)} \right] \nonumber \\
&={2 \over \pi} \cdot {\left( {1 \over 2}  \right)_a^3 \over 1_a^3} \cdot {1 \over 4^a} \cdot  \sum_{n=0}^{\infty} {1 \over 2^{2n}} \frac{ \left( a +\frac{1}{2}\right)_n^3}{(a+1)_n^3}[6(n+a)+1]. \nonumber
\end{align}
For proving (15) just use (14) and (2).
\end{proof}

\subsection{Identity 9}
Let $f(a)$ and $g(a)$ be the functions:
\[
f(a)=\sum_{n=0}^{\infty} {(-1)^{n} \over 2^{10n}} \frac{ \left( a + \frac{1}{2}\right)_n^5}{(a+1)_n^5} [820(n+a)^2+180(n+a)+13],
\]
\[
g(a)=\sum_{n=0}^{\infty} {1 \over 2^{6n}} \frac{ \left( a +\frac{1}{2}\right)_n^3}{(a+1)_n^3}[42(n+a)+5],
\]
then we have
\begin{equation}\label{}
f(a)=128a \sum_{n=0}^{\infty} \frac{ \left( a+\frac{1}{2} \right)_n^4}{(2a+1)_n^4} (4n+6a+1)
\end{equation}
and
\begin{equation}\label{}
f(a)= {8 \over \pi} \cdot {16^a \over \cos \pi a} \cdot { 1_a ^2 \over { \left(1 \over 2 \right)}_ a^2} \cdot g(a) +{2048 a^3 \over 2a-1} \sum_{n=0}^{\infty}  \frac{  \left( \frac{1}{2} +a \right)_n^3} {(2a+1)_n^2 \left( {3 \over 2} -a \right)_n}
\end{equation}
and also
\begin{align}
f(a) &={128 \over \pi^2} \cdot {1024^a \over \cos^3 \pi a} \cdot { 1_a ^5 \over { \left(1
\over 2 \right)}_ a^5} \\ &+ {1024 \over \pi} \cdot {16^a \over \cos \pi a} \cdot { 1_a ^2 \over { \left(1 \over 2 \right)}_ a^2} \cdot {a^2 \over 2a-1} \cdot \sum_{n=0}^{\infty}  \frac{  \left( a+\frac{1}{2} \right)_n^2} {(2a+1)_n \left( {3 \over 2} -a \right)_n}  \nonumber \\ &+{2048 a^3 \over 2a-1} \sum_{n=0}^{\infty}  \frac{ \left( \frac{1}{2} +a \right)_n^3} {(2a+1)_n^2 \left( {3 \over 2} -a \right)_n}. \nonumber
\end{align}
From (16) we easily get the evaluation [2]:
\[ f \left(\frac{1}{2}\right)=256 \zeta(3). \]
From (17) and using [1] the identity:
\[ \sum_{n=0}^{\infty} { {2n \choose n}^2 \over 16^n \cdot (2n+1)}= {4 G \over \pi}, \]
we get
\[
\lim_{a \to 0} \left[ {f(a) \over a^3}-{8 \over \pi} \cdot {16^a \over \cos \pi a} \cdot { 1_a ^2 \over { \left(1 \over 2 \right)}_ a^2} \cdot {g(a) \over a^3} \right]=-2048 \cdot \sum_{n=0}^{\infty} { {2n \choose n}^2 \over 16^n \cdot (2n+1)} =-8192 \cdot {G \over \pi}.
\]
And from (18) we get
\[
f(0)={128 \over \pi^2}, \quad \quad \quad f'(0)={2560 \over \pi^2} \ln 2, \quad \quad \quad
f''(0)={2560 \over 3 \pi^2} \left( 60 \ln^2 2 - \pi^2 \right).
\]
\begin{proof}
The proof of (15) is essentially encapsulated in the WZ-pair:
\[ F(n,k)=B(n,k) \cdot 128 \cdot n \cdot (6n+4k+1), \]
\begin{multline}
G(n,k)=B(n,k) \cdot \left[ {(2n+2k+1)^4 \over (2n+k+1)^4} \cdot (820n^2+180n+13+8k^2+20k+72nk) \right. \nonumber \\ -(296nk^3+1056n^2k^2+1280n^3k+528n^4+800n^3+1344n^2 k \nonumber \\ +608nk^2+28k^3+ \left. 408n^2+384nk+40k^2+72n+16k+1) \cdot {32nk \over (2n+k+1)^4} \right],
\end{multline}
where
\[
B(n,k)=\frac{(-1)^n}{2^{20n} \cdot 2^{8k}} \cdot \frac{(2n+2k)!^4 \cdot (2n)!^5}{(2n+k)!^4 \cdot (n+k)!^4 \cdot n!^6}
\]
and the proof of (16) is essentially encapsulated in the WZ-pair:
\[ F(n,k)=B(n,k) \cdot {2048 \cdot n^3 \over 2n-2k-1}, \]
\begin{multline}
G(n,k)=B(n,k) \cdot \left[ \frac{}{} 820n^2+180n+13+ {k \over (2n+k+1)^2}(1312n^3+1340n^2k+336nk^2 \right. \nonumber \\
\left. +1456n^2+828nk+40k^2+472n+79k+36) \frac{}{} \right], \nonumber
\end{multline}
where
\[
B(n,k)={(-1)^k \cdot (-1)^n \over 2^{20n} \cdot 2^{6k}} \cdot {(2n)!^4 \cdot (2n+2k)!^3
\cdot \left(n-k-{1 \over 2} \right)! \over n!^7 \cdot (n+k)!^3 \cdot (2n+k)!^2 \cdot \left(n-{1 \over 2} \right)!}.
\]
\end{proof}

\subsection{Identity 10}
Let $f(a)$ and $g(a)$ be the functions:
\[
f(a)=\sum_{n=0}^{\infty} {1 \over 2^{4n}} \frac{\left( a +\frac{1}{2}\right)_n^3 \left( a +\frac{1}{4}\right)_n
\left( a +\frac{3}{4}\right)_n}{(a+1)_n^5}[120(n+a)^2+34(n+a)+3)],
\]
\[
g(a)=\sum_{n=0}^{\infty} {1 \over 2^{6n}} \frac{ \left( a +\frac{1}{2}\right)_n^3}{(a+1)_n^3}[42(n+a)+5],
\]
then we have
\begin{equation}\label{}
f(a)=32a \sum_{n=0}^{\infty} \frac{ \left(\frac{1}{2} \right)_n^2 \left( a+\frac{1}{2} \right)_n^2}{(a+1)_ n^2 (2a+1)_n^2} (4n+4a+1)
\end{equation}
and
\begin{equation}\label{}
f(a)= {2 \over \pi} \cdot {1 \over 4^a \cdot \cos 2 \pi a} \cdot { 1_a ^2 \over { \left(1 \over 4 \right)}_ a { \left(3 \over 4 \right)}_ a} \cdot g(a) +{512 a^3 \over 4a-1} \sum_{n=0}^{\infty} \frac{  \left( \frac{1}{2} \right)_n^3} {(a+1)_n^2 \left( {3 \over 2} -2a \right)_n}
\end{equation}
and also
\begin{multline}
f(a)={32 \over \pi^2} \cdot {16^a \over \cos^2 \pi a } \cdot  {1 \over \cos 2 \pi a} \cdot { 1_a ^5 \over { \left(1
\over 2 \right)}_ a^3 \left(1 \over 4 \right)_ a \left(3 \over 4 \right)_ a }+{512 a^3 \over 4a-1} \sum_{n=0}^{\infty}  \frac{ \left( \frac{1}{2} \right)_n^3} {(a+1)_n^2 \left( {3 \over 2} -2a
\right)_n}
\\
+{256 \over \pi} \cdot {1 \over 4^a \cdot \cos 2 \pi a} \cdot {
1_a ^2 \over { \left(1 \over 4 \right)_ a \left(3 \over 4
\right)}_ a } \cdot {a^2 \over 2a-1} \cdot \sum_{n=0}^{\infty}
\frac{  \left( a+\frac{1}{2} \right)_n^2} {(2a+1)_n \left( {3
\over 2} -a \right)_n}.
\end{multline}
From (19) we easily get the evaluation:
\[ f \left(\frac{1}{2}\right)={16 \pi^2 \over 3}. \]
From (20) we get
\[
\lim_{a \to 0} \left[ {f(a) \over a^3} = {2 \over \pi} \cdot {1 \over 4^a \cdot \cos 2 \pi a} \cdot {
1_a ^2 \over { \left(1 \over 4 \right)}_ a { \left(3 \over 4 \right)}_ a} \cdot {g(a) \over a^3} \right]=-512 \cdot
\sum_{n=0}^{\infty} { {2n \choose n}^2 \over 16^n \cdot (2n+1)} =-2048 \cdot {G \over \pi}.
\]
And from (21) we get
\[ f(0)={32 \over \pi^2}, \quad \quad \quad f'(0)={512 \over \pi^2} \ln 2, \quad \quad \quad
f''(0)={64 \over 3 \pi^2} \left( 384 \ln^2 2 - 7 \pi^2 \right).
\]

\begin{proof}
The proof of (17) is essentially encapsulated in the WZ-pair:
\[ F(n,k)=B(n,k) \cdot 32n \cdot (4n+4k+1), \]
\begin{multline}
G(n,k)=B(n,k) \cdot \left[ {120n^2+34n+3+32k^2+128kn+16k \over 4n+4k+1} \nonumber \right. \\ \left. +k \cdot { 32n^3+8n^2k+16n+6kn+40n^2+k+2 \over (4n+4k+1)(2n+k+1)^2} \frac{}{} \right], \nonumber
\end{multline}
where
\[
B(n,k)=\frac{1}{2^{16n} \cdot 2^{8k}} \cdot \frac{(2k)!^2 \cdot (2n+2k)!^2 \cdot (4n)! \cdot (2n)!^2}
{(2n+k)!^2 \cdot (n+k)!^4 \cdot k!^2 \cdot n!^4}
\]
and the proof of (18) is essentially encapsulated in the WZ-pair:
\[ F(n,k)=B(n,k) \cdot {512 \cdot n^3 \over 4n-2k-1}, \]
\[ G(n,k)=B(n,k) \cdot (120n^2+84kn+34n+10k+3), \]
where
\[
B(n,k)={(-1)^k \over 2^{16n} \cdot 2^{6k}} \cdot {(2n)!^2 \cdot (4n)! \cdot
(2k)!^3 \cdot \left(2n-k-{1 \over 2} \right)! \over n!^6 \cdot
(n+k)!^2 \cdot k!^3 \cdot \left(2n-{1 \over 2} \right)!}.
\]
\end{proof}

\subsection*{Conclusion}
We have found a reduction of ten Ramanujan extended series to
simpler hypergeometric series. We believe all the other Ramanujan
extended series admit as well a reduction but we have not found
the adequate WZ-pairs to get them. It would also be interesting to
discover some reduction by using experimental techniques because
it could give the clue for finding the corresponding WZ-pairs
which would prove the identities.

\enddocument